\documentclass[12pt]{article}

\title{The inf-translation for solving set minimization problems}

\author{
Andreas H. Hamel\footnote{Free University of Bolzano-Bozen, \href{mailto:andreas.hamel@unibz.it}{andreas.hamel@unibz.it}}, Frank Heyde\footnote{Freiberg University of Mining and Technology, \href{mailto:frank.heyde@math.tu-freiberg.de}{frank.heyde@math.tu-freiberg.de}}, Daniela Visetti\footnote{Free University of Bolzano-Bozen, \href{mailto:daniela.visetti@unibz.it}{daniela.visetti@unibz.it}}
}
\date{{\small \today}}

\usepackage{array} 
\usepackage{verbatim} 

\usepackage{amsmath, amssymb, enumerate, color}
\usepackage{stmaryrd}

\newtheorem{theorem}{Theorem}
\newtheorem{corollary}[theorem]{Corollary}
\newtheorem{remark}[theorem]{Remark}
\newtheorem{lemma}[theorem]{Lemma}
\newtheorem{definition}[theorem]{Definition}
\newtheorem{proposition}[theorem]{Proposition}
\newtheorem{example}[theorem]{Example}

\numberwithin{equation}{section}  
\numberwithin{figure}{section}    
\numberwithin{table}{section}     
\numberwithin{theorem}{section}

\newcommand{\of}[1]{\ensuremath{\left( #1 \right)}}
\newcommand{\norm}[1]{\ensuremath{\left\| #1 \right\|}}

\newcommand{\cb}[1]{\ensuremath{ \left\{ #1 \right\} }}
\newcommand{\sqb}[1]{\ensuremath{ \left[ #1 \right] }}

\newcommand{\bs}{\backslash}

\newcommand{\pend}{ \hfill $\square$}

\newcommand{\vp}{\ensuremath{\varphi}}

\newcommand{\F}{\ensuremath{\mathcal{F}}}
\newcommand{\G}{\ensuremath{\mathcal{G}}}

\DeclareMathOperator*{\wMin}{wMin}

\newcommand{\R}{\mathrm{I\negthinspace R}}
\newcommand{\OLR}{\overline{\mathrm{I\negthinspace R}}}

\newcommand{\dom}{{\rm dom \,}}

\newcommand{\cl}{{\rm cl \,}}

\newcommand{\co}{{\rm co \,}}

\newcommand{\Int}{{\rm int\,}}

\newcommand{\Ima}{{\rm Im \,}}

\usepackage[
colorlinks=true, linkcolor=blue, citecolor=blue, urlcolor=blue, filecolor=blue
auskommentieren und unten pdfborder aktivieren]{hyperref}

\definecolor{color0}{gray}{.50}
\definecolor{color1}{rgb}{0,.2,.8}
\definecolor{color2}{rgb}{1,.2,0}
\definecolor{color3}{rgb}{.2,.7,.6}

\begin{document}
\maketitle
\begin{center}
{\em Dedicated to the memory of Prof. Hang-Chin Lai}
\end{center}

\begin{abstract}
Set- and vector-valued optimization problems can be re-formulated as complete lattice-valued problems. This has several advantages, one of which is the existence of a clear-cut solution concept which includes the attainment as the infimum (not present in traditional vector optimization theory) and minimality as two potentially different features. The task is to find a set which is large enough to generate the infimum but at the same time small enough to include only minimizers.

In this paper, optimality conditions for such sets based on the inf-translation are given within an abstract framework. The inf-translation reduces the solution set to a single point which in turn admits the application of more standard procedures. For functions with values in complete lattices of sets, scalarization results are provided where the focus is on convex problems. Vector optimization problems, in particular a vectorial calculus of variations problem, are discussed as examples.
\end{abstract}


\medskip
MSC Primary 49K27; Secondary 90C48, 06B23

\section{Complete lattice-valued optimization problems}

A complete lattice is a partially ordered set $(W, \leq)$ such that every subset of $W$ has an infimum and a supremum in $W$ with respect to $\leq$. Thus, the set $\R$ of real numbers together with the usual $\leq$ relation is not a complete lattice since, for example, $\R$ itself neither has an infimum nor a supremum in $\R$. On the other hand, $\OLR := \R \cup \{\pm\infty\}$ is a complete lattice where $\leq$ is extended in the obvious way. In fact, this complete lattice serves as a blueprint for many complete lattices of sets.

Lattices of sets sharing all properties with $(\overline \R, \leq)$ but the totalness of the order relation were recently studied in connection with optimization problems with a set-valued objective since they can be generated via so-called set relations. See \cite{Hamel05Habil, HamelEtAl15Incoll, KuroiwaTanakaTruong97NA, Loehne11Book}, for example. Moreover, vector optimization problems can be extended to complete lattice-valued problems in a straightforward manner and treated as  set optimization problems \cite{CrespiSchrage14ArX, HeydeLoehne11Opt, Loehne11Book}.

It turns out that the two features ``attaining a minimal value" and ``attaining the infimum" become two different concepts for complete lattice-valued minimization problems. This means, a point $\bar x \in X$ might be a minimizer of a complete lattice-valued function $f$, but the value $f(\bar x)$ is not the infimum of $f$, i.e., $f$ does not attain its infimum in $\bar x$--despite the fact that the infimum exists. Thus, looking for minimal function values and looking for function values which yield the infimum become two different tasks as already observed in \cite{Hamel04Con}.

In \cite[Definition 2.7]{HeydeLoehne11Opt}, \cite[Definition 2.8]{Loehne11Book}, solutions of complete lattice-valued optimization problems were defined which incorporate both minimality and the attainment of the infimum. The definition reads as follows.

\begin{definition}
\label{DefLatticeSolution}
Let $(W, \leq)$ be a complete lattice, $X \neq \emptyset$ a set and $f \colon X \to W$.

A set $M \subseteq X$ is called a {\em lattice-infimizer} of $f$ if
\[
\inf\cb{f(x) \mid x \in M} = \inf\cb{f(x) \mid x \in X}.
\]

An element $\bar x \in X$ is called a {\em lattice-minimizer} of $f$ if
\[
x \in X, \; f(x) \leq f(\bar x) \; \Rightarrow \; f(x) = f(\bar x). 
\]

A set $M \subseteq X$ is called a {\em lattice-solution} to the problem
\[
\tag{P} \text{minimize} \quad f \quad \text{over} \quad X
\]
if $M$ is a lattice-infimizer and each $x \in M$ is a lattice-minimizer of $f$. A lattice-solution $M \subseteq  X$ to (P) is called {\em full} if $M$ includes all lattice-minimizers of $f$.
\end{definition}

The labels $(W, \leq)$-infimizer and $(W, \leq)$-solution will also be used in the following if the lattice is emphasized in which an infimizer or a solution is looked for.

Note that a lattice-solution is called a mild solution in \cite[Definition 7.1]{HeydeLoehne11Opt} and \cite[Definition 2.41]{Loehne11Book}, and a full lattice-solution is just called a solution. Below, also a slightly different solution concept will be discussed which is the reason why the label ``(full) lattice-solution" is preferred in this paper. Moreover, the concept of minimizers also makes sense if $(W, \leq)$ is just a partially ordered set. In fact, it is the predominant solution concept in vector optimization (with respect to a vector order, also called efficient or minimal solution) and in many papers on set-valued optimization (with respect to a vector order or a set relation).

In a complete lattice, there always exist a top and a bottom element, sometimes dubbed $+\infty$ and $-\infty$, respectively. Thus, one can easily incorporate constraints in (P): if $\mathcal X \subseteq X$ is a set of feasible points, one can replace the values of $f$ by $+\infty$ outside $\mathcal X$. 

\begin{example}
\label{ExFandGLattice}
Typical examples of complete lattices relevant for set optimization problems can be constructed as follows. Let $Z$ be a non-trivial, topological linear space over the real numbers. Furthermore, let $C \subseteq Z$ be a non-empty convex cone (i.e., $C + C \subseteq C$ and $sC \subseteq C$ for all $s>0$) and
\[
\F(Z, C) = \cb{A \subseteq Z \mid A = \cl(A + C)}, \quad \G(Z, C) = \cb{A \subseteq Z \mid A = \cl\co(A + C)}
\]
where $\cl$ and $\co$ stands for the topological closure and the convex hull, respectively, and the addition of sets is the usual element-wise (Minkowski) addition with the extension $A + \emptyset = \emptyset + A = \emptyset$ for all $A \subseteq Z$. Then, both of $(\F(Z, C), \supseteq)$ and $(\G(Z, C), \supseteq)$ are complete lattices with top element $\emptyset$ and bottom element $Z$. If the cone is replaced by its negative, one obtains the complete lattices $(\F(Z, -C), \subseteq)$ and $(\G(Z, -C), \subseteq)$ with completely symmetric properties; in particular, $\emptyset$ is the bottom and $Z$ is the top element in these lattices.

Such complete lattices of sets can be generated as representatives of equivalence classes with respect to the symmetric part of so-called set relations (see the early work \cite{KuroiwaTanakaTruong97NA}, for example). Compare \cite{Hamel05Habil, HamelEtAl15Incoll} for the procedure and more and earlier references.
\end{example}

The solution concept in Definition \ref{DefLatticeSolution} may appear intricate since infimizers, in particular solutions of (P), are sets which are, in general, not singletons. Thus, optimality conditions such as ``zero belongs to the subdifferential of the function at some point" should actually be taken ``at sets" rather than ``at points" if one wants to characterize infimizers and solutions. 

It is already not straightforward to define the (directional) derivative or the subdifferential of a set-valued function at a point: several different concepts are used in vector optimization (especially with a set-valued objective) such as the contingent derivative \cite[Chap. 3]{DinhTheLuc89Book}, the contingent epiderivative \cite[Chap. 15]{Jahn04Book} and Mordukhovich's coderivatives \cite[Chap. 1]{Mordukhovich06Book1}, for example. They all share the feature that they are taken at points of the graph of a set-valued map rather than at points of the domain of a set-valued function. 

Clearly, it would be even more challenging to come up with such a derivative concept ``at a set" (instead of at a point of the domain or at points of the graph). A major goal of this paper is to provide and apply a general method to circumvent this difficulty in linear spaces. The first attempt in this direction is \cite{HamelSchrage14PJO}, where a directional derivative for convex set-valued functions is defined which is more in the spirit of traditional directional derivatives. It is also shown in \cite{HamelSchrage14PJO} that in the convex case this derivative based on residuations in complete lattices with an additional additive structure is not only a reformulation, but also a generalization of the coderivative from \cite[Chap. 1]{Mordukhovich06Book1}.

Another feature of the solution concept in Definition \ref{DefLatticeSolution} is emphasized by the obvious fact that $M \subseteq X$ is a lattice-infimizer of $f$ if, and only if,
\[
\tag{$\ell$-inf} N \; \text{is a lattice-infimizer of} \; f \; \text{for all} \; M \subseteq N \subseteq X.
\]
In particular, $X$ itself is always an infimizer of $f$. On the other hand, a single minimizer of $f$ often does not provide enough infomation as a solution of a vector- or set-valued optimization problem. Thus, one can turn the above solution concept into the task to find a set $M \subseteq X$ which is "big" enough to be an inifimizer, but also "small" enough to consist only of minimizers.

Under usual assumptions lattice-solutions of complete lattice-valued optimization problems exist. The standard reference is \cite{HeydeLoehne11Opt} (see also \cite{Loehne11Book}).  One more concept is needed to formulate the result.

\begin{definition}
\label{DefLevelClosedness} 
Let $X$ be a topological space and $(W, \leq)$ a partially ordered set. A function $f \colon X \to W$ is called level-closed if the set $L_f(w) = \cb{x \in X \mid f(x) \leq w}$ is closed for each $w \in W$.
\end{definition}

Level-closedness and a usual compactness assumption guarantee the existence of lattice-solutions.

\begin{theorem}
\label{ThmExistenceDomin}
Let $X \neq \emptyset$ be a compact topological space, $(W, \leq)$ a complete lattice and $f \colon X \to W$ a  level-closed function. Then problem (P) has a full lattice-solution.
\end{theorem}

{\sc Proof.} This follows from \cite[Proposition 2.8]{HeydeLoehne11Opt} and \cite[Proposition 5.15]{HeydeLoehne11Opt} in exactly the same way as \cite[Theorem 6.2]{HeydeLoehne11Opt} from these two propositions, and the details are therefore omitted. One may also compare \cite[Chap. 2]{Loehne11Book}. \pend

\medskip
Thus, this result and the method of its proof is essentially due to \cite{HeydeLoehne11Opt} although the final existence result Theorem \ref{ThmExistenceDomin} was not stated in \cite{HeydeLoehne11Opt, Loehne11Book} in this generality---only versions for vector-valued functions were given. An alternative approach based on scalarizing families can be found in \cite{CrespiHamelRoccaSchrage18ArX}.

In the next section, conlinear spaces are defined and a few crucial properties identified for subsequent use. Section \ref{SecInfTranslation} contains the central concept for reducing lattice-solutions to single points, namely the inf-translation. In Section \ref{SecSetLattices}, the most important lattices of sets are discussed, Section \ref{SecScalar} provides a scalarization procedure along with a corresponding solution concept and Section \ref{SecVectorOP} gives an application to vector optimization problems with a vectorial calculus of variations problem as an example. A few remarks on maximization conclude the paper.

\section{Ordered conlinear spaces}
\label{SecConlinSpaces}

The concept of an (ordered) conlinear space was introduced in \cite{Hamel05Habil}. It is motivated by the fact that it captures the algebraic structure which is preserved when passing from a(n) (ordered) linear space to its power set with the element-wise addition and multiplication by scalars. Compare \cite{HamelEtAl15Incoll} for a survey and more references.

\begin{definition}
\label{DefConlinearSpace} A nonempty set $W$ together with two algebraic operations $+
\colon W \times W \to W$ and $\cdot \colon \R_+ \times W \to W$ is called a conlinear
space provided that
\\
(C1) $(W, +)$ is a commutative semigroup with neutral element $\theta$,
\\
(C2) (i) $\forall w_1, w_2 \in W$, $\forall r \in \R_+$: $r \cdot (w_1 + w_2) = r\cdot w_1 + r \cdot w_2$, 

(ii) $\forall w \in W$, $\forall r, s \in \R_+$: $s \cdot (r\cdot w) = (sr) \cdot w$, 

(iii) $\forall w \in W${\em :} $1 \cdot w = w$, 

(iv) $\forall w \in W${\em :} $0 \cdot w = \theta$.
\end{definition}

\begin{definition}
\label{DefConlinearSubSpace} 
A non-empty subset $V \subseteq W$ of the conlinear space $(W, +, \cdot)$ is called a conlinear subspace of $W$ if 

(v) $v_1, v_2 \in V$ implies $v_1 + v_2 \in V$ and 

(vi) $v \in V$ and $t \geq 0$ imply $t \cdot v \in V$.
\end{definition}

If $V$ is a conlinear subspace of $(W, +, \cdot)$, then $(V, +, \cdot)$ is a conlinear space itself with the same zero element.

An element $w \in W$ of the conlinear space $(W, +, \cdot)$ is called convex if
\[
\forall s \in [0,1] \colon w = s w + (1-s)w.
\]
The zero element is always convex, and the set of convex elements in a conlinear space is a conlinear subspace \cite[Prop. 12]{Hamel05Habil}. 

\begin{definition}
\label{DefOrderedCS}
A conlinear space $(W, +, \cdot)$ together with a preorder $\preceq$ on $W$ is called preordered conlinear space provided
that 

(i) $w, w_1, w_2 \in W$, $w_1 \preceq w_2$ imply $w_1 + w \preceq w_2 + w$,

(ii) $w_1, w_2 \in W$, $w_1 \preceq w_2$, $r \in \R_+$ imply $r \cdot w_1 \preceq r\cdot w_2$.
\\
If $\preceq$ is additionally antisymmetric, $(W, +, \cdot, \preceq)$ is called ordered conlinear space.
\end{definition}

In certain ordered conlinear spaces all elements satisfy a ``one-sided" convexity condition which turns out to be an essential structural property. An ordered conlinear space is called semiconvex if the condition
\begin{equation}
\label{EqSemiconvex}
\forall w \in W, \forall s \in (0,1) \colon sw + (1-s) w \leq w
\end{equation}
is satisfied. Semiconcavity is defined symmetrically.

Another property links the complete-lattice structure and the addition. Let $\of{W, +, \cdot, \leq}$ be an ordered conlinear space such that $(W, \leq)$ is a complete lattice. Then, it is called inf-additive if
\begin{equation}
\label{EqInfAdd}
\inf \of{\{w\} + A} = w + \inf A
\end{equation}
for all $A \subseteq W$. Sup-additivity is defined similarly. It can be shown that inf-additivity is equivalent to several other properties, one of which is that the addition on $W$ admits an inf-residuation. Compare \cite{HamelEtAl15Incoll, HamelSchrage12} for more details.

\section{The inf-translation} 
\label{SecInfTranslation}

From now on, let $X$ be a non-trivial linear space over the real numbers. We consider problem (P) for the function $f \colon X \to W$ with values in the complete lattice $(W, \leq)$. The central concept of this paper is introduced in the following definition.

\begin{definition}
\label{DefInfTranslation} Let $M \subseteq X$ be a non-empty set. The function $\hat f(\cdot; M) \colon X \to W$ defined by
\begin{equation}
\label{EqInfTranslation}
\hat f(x; M) = \inf_{y \in M}f(x+y)
\end{equation}
is called the inf-translation of $f$ by $M$.
\end{definition}

The value $\hat f(x; M)$ coincides with the canonical extension of $f$ at $M + \{x\}$. The latter concept was introduced in \cite{HeydeLoehne11Opt, Loehne11Book}, but with a different goal in mind: the function $f$ was extended to a function on the power set of $X$. Here, the canonical extension taken at $M + \{x\}$ is considered as a function on $X$. This concept was introduced in \cite{HamelSchrage14PJO} for the special case of $\G(Z, C)$-valued functions.

Clearly, if $M = \{y\}$ is a singleton, then $\hat f(\cdot; M)$ is just a shift of $f$. However, if $M$ includes more than one element, even the inf-translation of a vector-valued function is a genuine set-valued one which is one more reason to discuss vector optimization problems within a complete lattice framework.

A few elementary, but important properties of the inf-translation are collected in the following lemma. Such a statement was already given in \cite[Lemma 5.8, Proposition 5.9]{HamelSchrage14PJO}, but only for (convex) $\G(Z, C)$-valued functions.

\begin{lemma}
\label{LemInfShiftProps}
For a set $M \subseteq X$, one has
\begin{itemize}

\item[(a)] if $M \subseteq N \subseteq X$ then $\hat f(x; N) \leq \hat f(x; M)$ for all $x \in X$,

\item[(b)] $\inf_{x \in X} f(x) = \inf_{x \in X} \hat f(x; M)$,

\item[(c)] the following features are equivalent:
\begin{itemize}
\item[(c1)]  $M$ is a lattice-infimizer of $f$,
\item[(c2)] $\{0\} \subseteq X$ is a lattice-infimizer of $\hat f(\cdot; M)$,
\item[(c3)] $\hat f(0; M) = \hat f(0; N)$ for all $M \subseteq N \subseteq X$, 
\item[(c4)] $\{0\} \subseteq X$ is a lattice-infimizer of $\hat f(\cdot; N)$ for all $M \subseteq N \subseteq X$.
\end{itemize}
\end{itemize}
\end{lemma}

Of course, ($\ell$-inf) gives another equivalent characterization for the list in (c). The following fact is an immediate consequence of Lemma \ref{LemInfShiftProps}.

\begin{corollary}
\label{CorConvexSetTrans} A set $M \subseteq X$ is a lattice-infimizer for $f$ if, and only if, $\cb{0}$ is a lattice-infimizer for $\hat f(\cdot; \co M)$ and $\hat f(0; M) = \hat f(0; \co M)$.
\end{corollary}

{\sc Proof of Lemma \ref{LemInfShiftProps}.} (a) is immediate from the definition of $\hat f(\cdot; M)$ as is (b) from
\[
\inf_{x \in X} f(x) = \inf_{x \in X}\inf_{y \in X} f(x + y) \leq \inf_{x \in X}\inf_{y \in M} f(x + y) = \inf_{x \in X} f(x).
\]

The equivalence of (c1) and (c2) follows from $\hat f(0; M) = \inf_{y \in M}f(y)$ and (b). (c1) implies (c3) since, if $M$ is a lattice-infimizer and $M \subseteq N$,
\[
\hat f(0; M) = \inf_{x \in M} f(x) \geq \inf_{x \in N} f(x)  \geq \inf_{x \in X} f(x) = \inf_{x \in M} f(x) = \hat f(0; M)
\]
by (a) where $\hat f(0; N) = \inf_{x \in N} f(x)$ by definition.

(c3) implies (c4) since, if $y \in X$ and $M \subseteq N$,
\[
\hat f(0; M) = \hat f(0; N) = \hat f(0; X) = \inf_{x \in X} f(x) \leq \inf_{x \in N}f(y + x) = \hat f(y; N).
\]

Finally, (c2) is immediate from (c4). \pend

In particular (c2) and (c4) of the previous lemma show that a lattice-infimizer (a set, in general) can be reduced to a single point via the inf-translation.

In order to deal with convex functions and sets, the framework needs to be specialized further. For the next result, let $W$ be a conlinear space which is partially ordered by $\leq$, i.e., $(W, +, \cdot, \leq)$ is an ordered conlinear space. It is still assumed that $(W, \leq)$ is a complete lattice. 

A function $f \colon X \to W$ is called convex if
\[
\forall  x, y \in X, \forall s \in (0,1) \colon f(sx + (1-s)y) \leq sf(x) + (1-s)f(y).
\]
The next lemma is a generalization of \cite[Lemma 5.8 $(\text{c})$]{HamelSchrage14PJO} which was stated for $\G(Z,C)$-valued functions.

\begin{lemma}
\label{LemInfShiftConvex} 
If $f$ and $M$ are convex, so is $\hat f(\cdot; M) \colon X \to W$.
\end{lemma}

{\sc Proof.} Take $s \in (0,1)$, $x_1, x_2 \in X$. One has $M = sM +
(1-s)M$ since $M$ is convex. This and the convexity of $f$ yield
\begin{align*}
\hat f(sx_1 + (1- s)x_2; M) & = \inf_{y \in M}f(sx_1+ (1- s)x_2 + y) \\
    & = \inf_{y_1, y_2 \in M}\sqb{f(s(x_1 + y_1) + (1- s)(x_2 + y_2))} \\
    & \leq \inf_{y_1, y_2 \in M}\sqb{sf(x_1 + y_1) + (1- s)f(x_2 + y_2)} \\
    & = s\hat f(x_1; M) + (1- s)\hat f(x_2; M).
\end{align*}
\pend

\begin{lemma}
\label{LemConvexImages} If $(W, +, \cdot, \leq)$ is semiconvex and $f \colon X \to W$ is a convex function, then 

(a) $f(x) \in W$ is a convex element for each $x \in X$,

(b) if, additionally, $(W, \leq)$ is inf-additive and $M \subseteq X$ is convex, then $\inf_{x \in M} f(x) \in W$ is a convex element.
\end{lemma}

{\sc Proof.} (a) By convexity of $f$, $f(x) = f(sx + (1-s)x) \leq sf(x) + (1-s)f(x)$ where \eqref{EqSemiconvex} yields the opposite inequality. Since $\leq$ is antisymmetric, equality follows which means that $f(x)$ is a convex element. 

(b) Set $I_f(M) = \inf_{x \in M} f(x) \in W$. By \eqref{EqSemiconvex}, one has $s I_f(M) + (1-s) I_f(M) \leq I_f(M)$ for all $s \in (0,1)$. Fix $s \in (0,1)$. Then
\[
\forall x, y \in M \colon I_f(M) \leq f(sx + (1-s)y) \leq sf(x) + (1-s)f(y).
\]
Taking the infimum on the right hand side over $x \in M$ and $y \in M$ consecutively one obtains $I_f(M) \leq s I_f(M) + (1-s) I_f(M)$ since one can use inf-additivity \eqref{EqInfAdd}. Antisymmetry of $\leq$ yields equality for $s \in (0,1)$ whereas the case $s \in \{0,1\}$ is immediate. \pend

Consequently, if \eqref{EqSemiconvex} is satisfied and $f$, $M$ are convex, then $\hat f(\cdot; M)$ has only convex values (as well as $f$).

Note that $\hat f\of{\cdot; M}$ is not convex in general even if $f$ is convex, and, of course, lattice-infimizers are not necessarily convex sets.

\section{Complete lattices of sets for set optimization}
\label{SecSetLattices}

Let $Z$ be a non-trivial, topological linear space over the real numbers. Furthermore, let $C \subseteq Z$ be a convex cone. As already mentioned in Example \ref{ExFandGLattice}, the pair $(\F(Z, C), \supseteq)$ is a complete lattice in which infima and suprema are given by
\begin{equation}
\label{EqFInfSup}
\inf_{A \in \mathcal A} A = \cl\bigcup_{A \in \mathcal A} A \quad \text{and} \quad \sup_{A \in \mathcal A} A = \bigcap_{A \in \mathcal A} A.
\end{equation}
The structure $(\F(Z, C), \oplus, \cdot)$ is a conlinear space where the addition is $A \oplus B = \cl(A + B)$ and the multiplication with non-negative numbers is defined element-wise except for $0 \cdot A = \cl C$ for all $A \in \F(Z, C)$, i.e., $\cl C$ is the zero element in $(\F(Z, C), \oplus, \cdot)$.

Adding the convex hull produces $\G(Z, C) = \cb{A \subseteq Z \mid A = \cl\co(A + C)}$ which is a conlinear subspace of $(\F(Z, C), \oplus, \cdot)$ and also a complete lattice with respect to $\supseteq$ with
\begin{equation}
\label{EqGInfSup}
\inf_{A \in \mathcal A} A = \cl\co\bigcup_{A \in \mathcal A} A \quad \text{and} \quad \sup_{A \in \mathcal A} A = \bigcap_{A \in \mathcal A} A
\end{equation}
for $\mathcal A \subseteq \G(Z, C)$. Since each element of $\G(Z, C)$ is a convex set, $\G(Z, C)$ only includes convex elements; in fact, it is precisely the conlinear subspace of convex elements of $\F(Z, C)$ (see Section \ref{SecConlinSpaces}).

Note that $(\F(Z, C), \supseteq)$ is semiconvex, but not semiconcave in general whereas the opposite is true for $(\F(Z, -C), \subseteq)$. Moreover, $(\G(Z, C), \supseteq)$ satisfies \eqref{EqSemiconvex} with equality instead of inequality. Finally, both $(\F(Z, C), \supseteq, \oplus)$ and $(\G(Z, C), \supseteq, \oplus)$ are inf-additive, but not sup-additive in general whereas the opposite is true for $(\F(Z, -C), \subseteq, \oplus)$ and $(\G(Z, -C), \subseteq, \oplus)$.

\begin{corollary}
\label{CorConvexFValued}
Let $f \colon X \to \F(Z, C)$ be a convex function and let $M \subseteq X$ be a convex set. Then

(a) $f(x) \in \G(Z, C)$ for all $x \in X$, 

(b) $\inf_{x \in M} f(x) \in \G(Z, C)$.
\end{corollary}

{\sc Proof.} Both claims follow from Lemma \ref{LemConvexImages} since $(\F(Z, C), \supseteq)$ satisfies \eqref{EqSemiconvex} and $(\F(Z, C), \oplus)$ is inf-additive. \pend

(a) means that, for a convex  $\F(Z, C)$-valued function, it is not a restriction to assume that it is $\G(Z, C)$-valued. (b) means that the infima of a convex function $f$ over convex sets (in particular, over $X$) in $(\F(Z, C), \supseteq)$ and $(\G(Z, C), \supseteq)$ coincide, and, as a consequence of \eqref{EqGInfSup},
\begin{equation}
\label{EqConvexInf}
\inf_{x \in X} f(x) =  \cl\bigcup_{x \in X} f(x) \in \G(Z, C).
\end{equation}
Of course, one can also verify directly that the set on the right hand side is convex.

\begin{remark}
\label{RemInfimizerFvsG}
If $f \colon X \to \F(Z, C)$ is a convex function and $M \subseteq X$ is an infimizer of $f$ in $(\F(Z, C), \supseteq)$, then $M$ also is an infimizer of $f$ in $(\G(Z, C), \supseteq)$ since in this case
\[
\inf_{x \in X} f(x) = \cl \bigcup_{x \in M} f(x) \subseteq  \cl\co \bigcup_{x \in M} f(x) \subseteq \inf_{x \in X} f(x)
\]
where the two outer infima are the same in $(\F(Z, C), \supseteq)$ and $(\G(Z, C), \supseteq)$.

On the other hand, if $M \subseteq X$ an infimizer of $f$ in $(\G(Z, C), \supseteq)$, then $\co M \subseteq X$ is an infimizer of $f$ in $(\F(Z, C), \supseteq)$ since the set
\[
\bigcup_{x \in \co M} f(x)
\]
is convex due to the convexity of $f$ and $\co M$ as one directly checks. Consequently, if $f$ and $M$ are convex, then $M$ is an infimizer of $f$ in $(\F(Z, C), \supseteq)$ if, and only if, it is an infimizer in $(\G(Z, C), \supseteq)$.
\end{remark}

\begin{remark}
\label{RemInfTransFvsG}
One should be equally careful with respect to the inf-translations. Since
\[
\cl\bigcup_{y \in M} f(x + y) \subseteq \cl\co\bigcup_{y \in M} f(x + y) 
\]
the inf-translation in $(\F(Z, C), \supseteq)$ is a subset of the inf-translation in $(\G(Z, C), \supseteq)$ (see \eqref{EqFInfSup}, \eqref{EqGInfSup} and the definition of the inf-translation \eqref{EqInfTranslation} in a general lattice). On the other hand, if both $f$ and $M$ are convex, then the two inf-translations coincide. The next simple example shows that the case of a convex $f$ and a non-convex $M$ deserves special attention.
\end{remark}

\begin{example}
Let $f \colon \R^2 \to \G(\R^2, \R^2_+)$ be defined by
\[
f(x) = \left\{
	\begin{array}{ccc}
	x + \R^2_+ & : & x \in \R^2_+, \; x_1 + x_2 \geq 1 \\
	\emptyset & : & \text{otherwise}
	\end{array}
	\right.
\]
One has
\[
\inf_{x \in \R^2} f(x) = \co\sqb{\of{\cb{(1,0)^T} + \R^2_+} \bigcup \of{\cb{(0,1)^T} + \R^2_+}},
\]
thus the finite set $M = \cb{(1,0)^T, (0,1)^T}$ is a $\G(\R^2, \R^2_+)$-infimizer of $f$. Of course, the features of this example are shared by functions $f$ with polyhedral graphs, e.g., objectives of linear vector optimization problems. Compare \cite[Chap. 4]{Loehne11Book} and \cite{Weissing20MMOR} for corresponding solution concepts.
\end{example}

\section{Scalarization results}
\label{SecScalar}

In this section, $Z$ is a non-trivial separated (Hausdorff) locally convex topological linear space over the real numbers. Its (non-trivial) topological dual is denoted $Z^*$. If $C \subseteq Z$ is a convex cone, then the set
\[
C^+ = \cb{z^* \in Z^* \mid \forall z \in C \colon z^*(z) \geq 0}
\]
denotes the (positive) dual cone of $C$. If there is an element $\hat z \in C$ such that $z^*(\hat z) > 0$ for all $z^* \in C^+\bs\{0\}$ (in particular, if $\Int C \neq \emptyset$), then the set
\[
B^+(\hat z) = \cb{y^* \in C^+ \mid y^*(\hat z) = 1}
\]
is a (closed and convex) base of $C^+$, i.e., for each element $z^* \in C^+\bs\{0\}$ there are unique $y^* \in B^+(\hat z)$ and $s > 0$ such that $z^* = s y^*$.

\subsection{The general procedure}

Let $f \colon X \to \F(Z, C)$ be a function. A family of extended real-valued functions $\vp_{f, z^*} \colon X \to \OLR$ with $z^* \in C^+$ is defined by
\begin{equation}
\label{EqScalarization}
\vp_{f, z^*}(x) = \inf_{z \in f(x)} z^*(z)
\end{equation}
where it is understood that $\vp_{f, z^*}(x) = +\infty$ whenever $f(x) = \emptyset$. Then, $\vp_{f, z^*}$ is convex for all $z^* \in C^+$ if, and only if, $f$ is convex. See \cite[Lemma 4.20 (a)]{HamelEtAl15Incoll} for a proof of this fact and earlier references. In results like this, $C^+$ can be replaced by a base $B^+(\hat z)$ of $C^+$ if it exists since the functions $\vp_{f, z^*}$ are positive homogeneous in $z^*$.

Note that $\vp_{f, z^*}(x) = \inf_{z \in \cl\co f(x)} z^*(z)$, thus the family $\cb{\vp_{f, z^*}}_{z^* \in C^+}$ does not distinguish between a function $f \colon X \to \F(Z, C)$ and the $\G(Z, C)$-valued function for which every value $f(x)$ is replaced by its closed convex hull. This is due to the fact that $\vp_{f, z^*}$ is a version of the support function of the function values of $f$ and hints that this type of scalarization is appropriate for convex $f$. However, note that the infima in $\F(Z, C)$ and $\G(Z, C)$ can be different as \eqref{EqFInfSup} and \eqref{EqGInfSup} show.

\begin{lemma}
\label{LemInfimizerScalar}
(a) Let $f \colon X \to \F(Z, C)$ be a function. If $M \subseteq X$ is a lattice-infimizer of $f$, then $0 \in X$ is a minimizer of $\vp_{\hat f(\cdot; M), z^*}$ for every $z^* \in C^+\bs\{0\}$. The converse is true if $f$ and $M$ are convex.

(b) Let $f \colon X \to \G(Z, C)$ be a function. If $M \subseteq X$ is a lattice-infimizer of $f$, then $0 \in X$ is a minimizer of $\vp_{\hat f(\cdot; M), z^*}$ for every $z^* \in C^+\bs\{0\}$. The converse is true if $f$ is convex.
\end{lemma}

{\sc Proof.} (a) First, let $M \subseteq X$ be an infimizer of $f$ in $(\F(Z,C), \supseteq)$. Take $z^* \in C^+\bs\{0\}$ and $x \in X$. Then
\begin{align*}
\vp_{\hat f(\cdot; M), z^*}(x) & =  \inf\cb{z^*(z) \mid z \in \hat f(x; M)} = \inf\cb{z^*(z) \mid z \in \inf_{y\in M}f(x + y)} \\
	& \geq \inf\cb{z^*(z) \mid z \in \inf_{y \in X}f(y)} = \inf\cb{z^*(z) \mid z \in \inf_{y \in M}f(y)} \\
	& =  \inf\cb{z^*(z) \mid z \in \hat f(0, M)} = \vp_{\hat f(\cdot; M), z^*}(0)
\end{align*}
hence $0$ is a minimizer of $\vp_{\hat f(\cdot; M), z^*}$.

For the converse, assume $f$ and $M$ are convex and that $0 \in X$ is a minimizer of $\vp_{\hat f(\cdot; M), z^*}$ for every $z^* \in C^+\bs\{0\}$. Suppose $M$ is not an infimizer of $f$. Then, there exist $\bar x \in X$ and $\bar z \in f(\bar x)$ such that 
\[
\bar z \not \in \inf_{y \in M} f(y) = \hat f(0; M)
\]
since otherwise
\[
\bigcup_{x \in X} f(x) \subseteq \inf_{y \in M} f(y) = \hat f(0; M)
\]
and hence
\[
\inf_{x \in X} f(x) = \cl\of{\bigcup_{x \in X} f(x)} \subseteq \inf_{y \in M} f(y)
\]
which means that $M$ would be an infimizer. Since $\hat f(0; M) = \cl\bigcup_{y \in M}f(y)$ is a closed convex set due to the convexity of $M$ and $f$ (see Remark \ref{RemInfimizerFvsG}), $\bar z$ can be strongly separated from it, i.e., there is $z^* \in Z^*\bs\{0\}$ satisfying
\[
z^*(\bar z) < \inf_{z \in \hat f(0; M)} z^*(z).
\]
Since $f$ and $M$ are convex, so is $\hat f(\cdot; M)$ with values in $\G(Z, C)$. Therefore, one has $z^* \in C^+\bs\{0\}$. Pick $\bar y \in M$. Then
\begin{align*}
\vp_{\hat f(\cdot; M), z^*}(\bar x - \bar y) & = \inf\{z^*(z) \mid z \in \hat f(\bar x - \bar y; M)\} \\
	& \leq  \inf\{z^*(z) \mid z \in f(\bar x)\} \leq z^*(\bar z) \\
	& < \inf_{z \in \hat f(0; M)}z^*(z) = \vp_{\hat f(\cdot; M), z^*}(0).
\end{align*}
This contradicts the assumption that $0 \in X$ is a minimizer of $\vp_{\hat f(\cdot; M), z^*}$. 

(b) For the first claim, observe that the same estimate works for the infimum and the inf-translation in $(\G(Z,C), \supseteq)$.

For the converse direction, observe that now $\hat f(0; M) = \cl\co\bigcup_{y \in M}f(y)$ is a closed convex set by definition, so the same separation argument as in (a) works. \pend

Lemma \ref{LemInfimizerScalar} can be understood as a necessary condition for infimizers. If convexity is present, the condition becomes also sufficient. The condition that $0$ is a minimizer of $\vp_{\hat f(\cdot; M)}$ can be directly characterized via $\vp_{f, z^*}$.

Next, the relationship between the scalarization $\vp_{\hat f(\cdot; M), z^*}$ of the inf-translation of $f$ and the inf-translation $\hat \vp_{f, z^*}(\cdot; M) \colon X \to \OLR$ of $\vp_{f, z^*}$ defined by
\[
\hat \vp_{f, z^*}(x; M) = \inf_{y \in M}\vp_{f, z^*}(x + y)
\]
is clarified. In fact, the two operations ``taking the inf-translation" and ``scalarization via \eqref{EqScalarization}" commute.

\begin{proposition}
\label{PropTransScalar}
Let $f \colon X \to \F(Z, C)$ or $f \colon X \to \G(Z, C)$ be a function. Then, for all $z^* \in C^+$ one has
\[
\forall x \in X \colon \vp_{\hat f(\cdot; M), z^*}(x) = \hat \vp_{f, z^*}(x; M).
\]
\end{proposition}

{\sc Proof.}  The definitions and the continuity of $z^* \in C^+$ yield
\begin{align*}
\vp_{\hat f(\cdot; M), z^*}(x) & = \inf\cb{z^*(z) \mid z \in \hat f(x; M)} =  \inf\cb{z^*(z) \mid z \in \cl\bigcup_{y \in M} f(x + y)} \\
	& = \inf_{y \in M} \inf\cb{z^*(z) \mid z \in f(x + y)} = \inf_{y \in M} \vp_{f, z^*}(x + y) = \hat \vp_{f, z^*}(x; M),
\end{align*}
and the same lines work if $\cl\bigcup_{y \in M} f(x + y)$ is replaced by $\cl\co\bigcup_{y \in M} f(x + y)$. This completes the proof. \pend

\begin{lemma}
\label{LemScalarInfimizer} Let $f \colon X \to \F(Z, C)$ be a function, $M \subseteq X$ and $z^* \in C^+\bs\{0\}$. Then 
\begin{equation}
\label{EqScalarInfimizer}
\inf_{y \in M} \vp_{f, z^*}(y) = \inf_{x \in X}\vp_{f, z^*}(x)
\end{equation}
if, and only if, $0 \in X$ is a minimizer of $\vp_{\hat f(\cdot; M), z^*}$.
\end{lemma}

{\sc Proof.} This follows from 
\[
\vp_{\hat f(\cdot; M), z^*}(0) = \inf_{y \in M}  \vp_{f, z^*}(y) = \inf_{x \in X}  \vp_{f, z^*}(x) = 
	\inf_{x \in X} \inf_{y \in M} \vp_{f, z^*}(x + y) =  \inf_{x \in X} \vp_{\hat f(\cdot; M), z^*}(x)
\]
and Proposition \ref{PropTransScalar} since the second term is the definition of $\hat \vp_{f, z^*}(0; M)$ which is equal to  $\vp_{\hat f(\cdot; M), z^*}(0)$ by this proposition. \pend

Clearly, condition \eqref{EqScalarInfimizer} is nothing but the fact that $M$ is an infimizer for $\vp_{f, z^*}$. Lemma \ref{LemScalarInfimizer} and Proposition \ref{PropTransScalar} together just mean that $0$ is a minimizer of $\hat \vp_{f, z^*}(x; M)$ if, and only if, $M$ is an infimizer of $\vp_{f, z^*}$---which of course also is a consequence of Lemma \ref{LemInfShiftProps} applied to $\vp_{f, z^*}$ and $\hat \vp_{f, z^*}(\cdot; M) = \vp_{\hat f(\cdot; M), z^*}$. 

\begin{corollary}
\label{CorScinfimizer}
(a) Let $f \colon X \to \F(Z, C)$ and $M \subseteq X$ be convex. Then $M$ is a lattice-infimizer of $f$ if, and only if, \eqref{EqScalarInfimizer} holds.

(b) Let $f \colon X \to \G(Z, C)$ be convex and $M \subseteq X$. Then $M$ is a lattice-infimizer of $f$ if, and only if, \eqref{EqScalarInfimizer} holds.
\end{corollary}

{\sc Proof.} Both versions follow from Lemma \ref{LemScalarInfimizer} and Lemma \ref{LemInfimizerScalar} \pend

\begin{corollary}
\label{CorTransScalarConvex}
If $f$ and $M$ are convex, then $\hat \vp_{f, z^*}(\cdot; M) \colon X \to \OLR$ is convex.
\end{corollary}

{\sc Proof.} This follows from Proposition \ref{PropTransScalar}, Lemma \ref{LemInfShiftConvex} (a) and the fact that \eqref{EqScalarization} produces a convex function if $f$ is convex. \pend

Taking Proposition \ref{PropTransScalar} into account, one can use Lemma \ref{LemInfimizerScalar} via first order conditions for minimizing $\hat \vp_{f, z^*}(x; M)$ over $x \in X$: find a set $M \subseteq X$ such that $\hat \vp_{f, z^*}(\cdot; M)$ satisfies a first order condition at $0 \in X$ for all $z^* \in C^+\bs\{0\}$.

Note that 
\[
\hat \vp_{f, z^*}(x; X) = \inf_{y \in X}\vp_{f, z^*}(x + y) \equiv \inf_{y \in X}\vp_{f, z^*}(y)
\]
is a constant function, thus every point is a minimizer. Thus, Lemma \ref{EqScalarization} becomes trivial---as $M=X$ always is a lattice-infimizer of $f$. In light of Definition \ref{DefLatticeSolution}, the task is to find an infimizer which is as ``small" as possible, e.g., only includes lattice-minimizers. The next subsection presents another possibility which has a close link to weakly minimal solutions in vector optimization.

\subsection{Solutions of scalarized problems and optimality conditions}

The next concept is concerned with a stronger version of condition \eqref{EqScalarInfimizer}.

\begin{definition}
\label{DefScalarSolutions}
Let $f \colon X \to \F(Z, C)$ be a function and $z^* \in C^+\bs\{0\}$. A point $\bar x \in X$ is called a $z^*$-minimizer of $f$ if 
\[
\forall x \in X \colon \vp_{f, z^*}(\bar x) \leq \vp_{f, z^*}(x).
\]
A set $M \subseteq X$ is called a sc-solution (short for scalarization-solution) of (P) if it is a lattice-infimizer and only includes $z^*$-minimizers.
\end{definition}

This solution concept has been considered first in \cite{HamelSchrage14PJO}. Condition \eqref{EqScalarInfimizer} is satisfied if the function $\vp_{f, z^*}$ has a minimizer for every $z^* \in C^+\bs\{0\}$. The assumption that such $z^*$-minimizers exists for every $z^* \in C^+\bs\{0\}$ is certainly a strong one. However, sc-solutions can be smaller.

\begin{example}
\label{ExHyperbolic}
Let $Z = \R^2$, $C=\R^2_+$ and $f \colon \R \to \G(Z,C)$ be defined by
\[
f(x) =
\left\{
	\begin{array}{ccc}
		\of{x, \frac{1}{x}}^T + \R^2_+ & : &  x > 0\\
		\emptyset & : & x \leq 0
	\end{array}
\right.	
\]
Neither for $z^* = (1,0)^T$, nor for $z^* = (0,1)^T$ does a $z^*$-minimizer exist. On the other hand, the set $M = \{x \in \R \mid x > 0\}$ is a sc-solution. This shows that the existence of a solution of every scalarized problem is not necessary.
\end{example}

We state an optimality condition for convex problems: such a condition is necessary and sufficient due to convexity. However, other versions based on different types of subdifferentials or derivatives could also be formulated. Then, the not-so-easy task is to ensure appropriate conditions to functions of the type $x \mapsto \vp_{f, z^*}(x)$, $x \mapsto \hat \vp_{f, z^*}(x; M)$ such as Lipschitz continuity or even differentiability. Yet another task is to compute the corresponding subdifferentials and derivatives, respectively.

The subdifferentials in the following result are the usual ones in convex analysis. A standard reference is \cite{Zalinescu02Book}. One should note that up to now, $X$ was merely assumed to be a linear space. Stronger assumptions are only needed for optimality conditions.

\begin{theorem}[optimality condition]
\label{ThmOC}
Let $X$ be a separated, locally convex, topological linear space, let $f \colon X \to \G(Z, C)$ be a convex function and let $M \subseteq X$ be a set. Then, $M$ is a sc-solution of (P) if, and only if,

(1) $\inf_{y \in M}f(y) = \inf_{y \in \co M}f(y)$,

(2) one has
\[
\forall z^* \in C^+\bs\{0\} \colon 0 \in \partial \hat \vp_{f, z^*}(0; \co M),
\]

(3) for each $y \in M$ there is a $z^* \in C^+\bs\{0\}$ such that 
\[
0 \in \partial\vp_{f, z^*}(y).
\]
\end{theorem}

{\sc Proof.} By Corollary \ref{CorConvexSetTrans} and Lemma \ref{LemInfimizerScalar}, $M$ is a lattice-infimizer of $f$ if, and only if, the conditions in (1), (2) are satisfied since the functions $\hat \vp_{f, z^*}(\cdot; \co M)$ are convex. The convexity of the functions $\vp_{f, z^*}$ ensures that (3) holds if, and only if, each $y \in M$ is a $z^*$-minimizer. \pend

\begin{remark}
\label{RemOneDim} If $Z = \R$, $C = \R_+$, then $z^* = 1$ is the only interesting element since $\{1\}$ is a base of $C^+ = \R_+$. Set $\vp(x) = \inf f(x)$. Moreover, attainment of the infimum at a point $\bar x \in X$ and minimality of $\bar x \in X$ are equivalent, and so are the two conditions in (2), (3) for $M = \{\bar x\}$ since $\vp_{f, 1}(\bar x) = \vp(\bar x)$ and
\[
\hat \vp_{f, 1}(0; \{\bar x\}) = \vp_{f, 1}(\bar x) = \vp(\bar x).
\]
\end{remark}

The previous result has a similar algorithmic character as a scalar necessary optimality condition such as Fermat's rule: find $z^*$-minimizer for as many $z^* \in C^+\bs\{0\}$ as possible (this is a parametric scalar optimization problem) and then check if the set of these $z^*$-minimizers already is a lattice-infimizer, i.e., if the conditions (1), (2) of Theorem \ref{ThmOC} are satisfied. As in the scalar case, necessary conditions become sufficient if convexity is present.

\begin{example}[Example \ref{ExHyperbolic} cont.]
\label{ExHyperbolicCont}
It is sufficient to determine $\hat \vp_{f, w_\alpha}(x; M)$ for $w_\alpha = \alpha(1,0)^T + (1-\alpha)(0,1)^T \in \R^2_+ = C^+$ with $\alpha \in [0,1]$ since $\R^2_+$ has a base generated by $\bar z = (1, 1)^T$. The set $M = \{y \in \R \mid y > 0\}$ (see  Example \ref{ExHyperbolic}) is the set of all $w_\alpha$-minimizers for $\alpha \in (0,1)$ and an infimizer according to Theorem \ref{ThmOC}. One can compute
\[
\hat \vp_{f, w_0}(x; M) \equiv 0, \quad \hat \vp_{f, w_1}(x; M) = x
\]
for all $x > 0$ and
\begin{align*}
\hat \vp_{f, w_\alpha}(x; M) & = \inf_{y > 0} \sqb{\alpha(x+y) + (1-\alpha)\frac{1}{x+y}} \\
	& = \left\{
		\begin{array}{ccc}
		2\sqrt{\alpha(1-\alpha)} & : & x \leq \sqrt{\frac{1-\alpha}{\alpha}} \\
		\alpha x + (1-\alpha)\frac{1}{x} & : & x \geq \sqrt{\frac{1-\alpha}{\alpha}}
		\end{array}
		\right.
\end{align*}
for $\alpha \in (0,1)$ and for all $x > 0$. Moreover,
\[
\forall \alpha \in [0,1] \colon \hat \vp_{f, w_\alpha}(0; M) = 2\sqrt{\alpha(1-\alpha)}.
\]
One can easily check that $\hat \vp_{f, w_\alpha}(x; M)$ is increasing and convex for $x \geq 0$, hence $\{0\} \subset \R$ is an infimizer of $\hat \vp_{f, w_\alpha}(\cdot; M)$; moreover, $0 \in \partial \hat \vp_{f, w_\alpha}(0; M)$ in all cases in accordance with Theorem \ref{ThmOC}. Note that $\vp_{f, w}(0) = +\infty$ since $f(0) = \emptyset$, but $\hat \vp_{f, w_\alpha}(0; M) \in \R$.
\end{example}

The following result provides a different type of optimality condition which will be used in the subsequent section.

\begin{proposition}\label{PropDiffGen}
Let $(X, \norm{\cdot})$ be a normed space, $f \colon X \to \F(Z, C)$, $M \subseteq X$  and $z^*\in C^+\bs\{0\}$. If there is $\bar x \in M$ and a linear subspace $D$ of $X$ with $M\subseteq\cb{\bar x} + D$ such that
\[
\vp_{f, z^*}(\bar x) = \min_{u \in D} \vp_{f, z^*}(\bar x+ u)
\]
and if $\psi \colon D\to\R$, defined by $\psi(u)=\vp_{f, z^*}(\bar x +u)$, is Fr{\'e}chet-differentiable at $0$, then $\hat \vp_{f, z^*}(\cdot; M)$ is Fr{\'e}chet-differentiable at $0$ on $D$ and one has 
\[
\hat \vp_{f, z^*}'(0; M)(u)=\psi'(0)(u)=0
\] 
for all $u\in D$.
\end{proposition}

{\sc Proof.} The assumptions imply
$M+D=\cb{\bar x}+D$, hence
\begin{align*}
  \inf_{y \in M}\vp_{f, z^*}(y)\le \vp_{f, z^*}(\bar x)=\inf_{u \in D}\vp_{f, z^*}(\bar x + u)= \inf_{u \in D}\inf_{y \in M}\vp_{f, z^*}(y + u) \le \inf_{y \in M}\vp_{f, z^*}(y).
\end{align*}
Therefore, one has
\begin{align*}
0 & \leq \hat \vp_{f, z^*}(u; M) - \hat \vp_{f, z^*}(0; M) = \inf_{y\in M}\vp_{f, z^*}(u+y) - \inf_{y \in M}\vp_{f, z^*}(y)\\
& = \inf_{y\in M}\vp_{f, z^*}(u+y) - \vp_{f, z^*}(\bar x) \leq \vp_{f, z^*}(u+\bar x) - \vp_{f, z^*}(\bar x)\\
& = \psi'(0)(u) + o(\norm{u}) = o(\norm{u})
\end{align*}
for all $u \in D$ since $0$ minimizes $\psi$ on $D$ and hence $\psi'(0)(u)=0$ for all $u \in D$. The statement follows. \pend

Again, with additional convexity assumptions, Lemma \ref{LemInfimizerScalar} turns the previous proposition into a sufficient condition for sc-solutions.

\section{Vector optimization problems}
\label{SecVectorOP}

In this section, let $X$ and $Z$ be (non-trivial) separated locally convex, topological vector spaces over the reals and $C \subseteq Z$ a closed convex pointed cone. As before, $X^*$, $Z^*$ denote the topological duals.

For a function $F \colon X \to Z \cup\{+\infty\}$, let the two sets
\[
\dom F = \cb{x \in X \mid F(x) \neq +\infty} \quad \text{and} \quad \Ima F =  \cb{F(x) \mid x \in \dom F}
\]
be the domain and the image of $F$. Such a function $F$ is called $C$-convex if
\[
s \in (0,1), \; x_1, x_2 \in \dom F \; \Rightarrow \; sF(x_1) + (1-s)F(x_2) \in \{F(sx_1 + (1-s)x_2)\} + C.
\]
A point $\bar x \in \dom F$ is called a minimizer of $F$ if
\[
\of{F(\bar x) - C} \cap \Ima F = \cb{F(\bar x)},
\]
and it is called a weak minimizer if
\[
\of{F(\bar x) - \Int C} \cap \Ima F = \emptyset
\]
where $\Int C \neq \emptyset $ is assumed. The set of weak minimizers of $F$ is denoted by $\wMin F$. This solution concept is very popular in vector optimization (see, for example, \cite[Definition 2.1 in Chap. 2]{DinhTheLuc89Book}), and weak minimizers are often called weakly efficient solutions. A basic fact links weak minimizers to solutions of scalarized problems.

\begin{lemma}
\label{LemScalarWeakMin}
Let $F$ be $C$-convex. Then $\bar x \in \wMin F$ if, and only if, there is $z^* \in C^+\bs\{0\}$ such that $\bar x$ is a minimizer of the function
$z^* \circ F \colon X \to \R\cup\{+\infty\}$ defined by
\[
(z^* \circ F)(x) = z^*(F(x))
\]
with the convention $(z^* \circ F)(x) = +\infty$ for $x \not\in \dom F$.
\end{lemma}

{\sc Proof.} This is \cite[Theorem 2.10 in Chap. 4]{DinhTheLuc89Book}. \pend

The link to set optimization is provided by the next concept. 

\begin{definition}
\label{DefSetExtension}
The inf-extension of a function $F \colon X \to Z \cup\{+\infty\}$ is the function $f \colon X \to \G(Z, C)$ defined by
\[
f(x) = \left\{
		\begin{array}{ccc}
		\{F(x)\} \oplus C  & : & x \in \dom F \\
		\emptyset  & : & x\not\in \dom F
		\end{array}
		\right.
\]
\end{definition}

It can easily be shown that $f$ is convex if, and only if, $F$ is $C$-convex. Moreover, minimizers of $F$ with respect to $\leq_C$ are one-to-one with minimizers of its inf-extension $f$ with respect to $\supseteq$ (see \cite{HeydeLoehne11Opt, Loehne11Book}).

\begin{definition}
\label{DefVectorSol}
A nonempty set $M \subseteq X$ is called a lattice-infimizer (lattice-solution, sc-solution) of the vector optimization problem
\[
\tag{VOP} \text{minimize} \quad F \quad \text{over} \quad X \quad \text{with respect to} \quad \leq_C 
\]
if it is a lattice-infimizer (lattice-solution, sc-solution) of the set optimization problem (P) where $f$ is the inf-extension of $F$.
\end{definition}

Compare \cite{CrespiSchrage14ArX} for such a set optimization approach to vector optimization problems. Note that these solution concepts are considered in $(\G(Z,C), \supseteq)$ with the corresponding formula for the infimum. Of course, one may define other types of lattice-infimizers and -solutions for (VOP). Since the focus in this paper is on $z^*$-minimizers which make mostly sense in a convex framework, alternatives are not discussed here.

\begin{corollary}
\label{CorVectorScalar} 
Let $F$ be $C$-convex. Then 

(a) $\bar x \in \wMin F$ if, and only if, $\bar x$ is a $z^*$-minimizer for the inf-extension $f$ of $F$ for some $z^* \in C^+\bs\{0\}$.

(b) if (VOP) has a $z^*$-minimizer for each $z^* \in C^+\bs\{0\}$, then the set $\wMin F$ is a sc-solution of (VOP).
\end{corollary}

{\sc Proof.} (a) This follows from Lemma \ref{LemScalarWeakMin} and 
\[
\vp_{f, z^*}(x) = (z^* \circ F)(x) = z^*(F(x)).
\]

(b) This follows from Corollary \ref{CorScinfimizer} with Lemma \ref{LemScalarWeakMin} in view.
\pend

Compare \cite[Proposition 2.15]{CrespiSchrage14ArX} for a related result. Two comments are in order. 

First, even though $z^*$-minimizers of a $C$-convex (VOP) with $\Int C \neq \emptyset$ are nothing but the well-known weakly efficient solutions, the inf-translation $\hat f(\cdot, M)$ is a "true" set-valued function in general---in particular for $M = \wMin F$. 

Secondly, the assumption that there is a $z^*$-minimizer for each $z^* \in C^+\bs\{0\}$ is a strong one as Example \ref{ExHyperbolic} shows. Therefore, it makes sense to look for a set $M \subseteq \wMin F$ which already produces the infimum. This brings---again, the "truly" set-valued---condition (2) of Theorem \ref{ThmOC} into play. For the inf-translation of $\vp_{f, z^*}$ one has in case of (VOP)
\[
\hat \vp_{f, z^*}(x; M) = \inf_{y \in M} \vp_{f, z^*}(x + y) = \inf_{y \in M} z^*(F(x + y)).
\]

\begin{example}[Example \ref{ExHyperbolic}, \ref{ExHyperbolicCont} cont.]
The function $\hat f(\cdot, M)$ for $M = \{y \in \R \mid y > 0\}$ in Example \ref{ExHyperbolic} has upper level sets of hyperbolas as values since
\[
\hat f(x, M) = \inf\cb{\cb{\of{x+y, \frac{1}{x+y}}^T} + \R^2_+ \mid y > 0} = \cl\bigcup_{y >0} \cb{\of{x+y, \frac{1}{x+y}}^T} + \R^2_+
\]
for $x \geq 0$. Example \ref{ExHyperbolicCont} provides the formula for $\hat \vp_{f, z^*}(x; M)$. Of course, $f$ is the inf-extension of $F \colon \R \to \R^2$ defined by
\[
F(x) =
\left\{
	\begin{array}{ccc}
		\of{x, \frac{1}{x}}^T & : &  x > 0\\
		+\infty & : & x \leq 0
	\end{array}
\right.	
\]
\end{example}

As a more complex example, a vector-valued calculus of variations problem is considered. Such problems arise in applications even in areas like the design of energy-saving buildings in architecture \cite{Marks97BE}. A more elaborate treatment of such problems with set optimization methods can be found in \cite{HamelVisetti18JMAA, HeydeVisetti19ArX}.

Let $a, b \in \R$ be two real numbers with $a < b$, $n, d$ two positive integers and $L \colon [a,b] \times \R^n \times \R^n \to \R^d$ be a function of class $C^1$.  We consider
\begin{equation}
\label{EqObjective}
F(x) = \int_a^b L(t,x(t),\dot x(t))\, dt
\end{equation}
on the set of feasible arcs
\begin{equation}
\label{EqConstraint}
\mathcal{X}= \{x \in C^1([a,b]; \R^n) \mid x(a)=A,\, x(b)=B\}
\end{equation}
for $A, B \in \R^n$. Since $L$ maps into $\R^d$ where $d$ can be strictly greater than 1, the problem of minimizing $F$ is a multi-criteria calculus of variations problem. 

It is assumed in the following that the function $L(t, \cdot, \cdot)$ is (jointly) convex in $(x, \dot x)$ for all $t \in [a,b]$. The set-valued extension of the problem requires to extend $F$ to a function mapping into $\G(\R^d, C)$. We define $f \colon C^1([a,b]; \R^n) \to \G(\R^d,C)$ by
\[
f(x) = \left\{
		\begin{array}{ccc}
		\{F(x)\} + C  & : & x \in \mathcal X \\
		\emptyset  & : & x\not\in \mathcal X
		\end{array}
		\right.
\]
The convexity assumption ensures that $F$ is $C$-convex and hence $f$ is convex. Now, the set-valued calculus of variations problem
\[
\tag{CVP} \text{minimize} \quad f(x) \quad \text{subject to} \quad x \in \mathcal X
\]
is considered. The inf-translation of $f$ with $M \subseteq C^1([a,b]; \R^n)$ is
\begin{equation}
\label{EqInfTransCVP}
\hat f(x; M) = \inf_{u \in M} f(x + u).
\end{equation}
One has $\hat f(x; M) = \emptyset$ whenever $M \subseteq \mathcal X$ and $x \not\in C^1_0([a,b]; \R^n)$ since in this case $x + u \not\in \mathcal X$ for all $u \in M$. This means that the inf-translation can be considered as a function on the linear subspace $C^1_0([a,b]; \R^n) $ of $C^1([a,b]; \R^n)$ which ``absorbs" the boundary conditions.

Clearly, as before, even though the original problem is vector-valued, the problem
\[
\tag{IP} \text{minimize} \quad \hat f(x; M) \quad \text{subject to} \quad x \in C^1_0([a,b]; \R^n)
\]
is a genuine set optimization problem, and one has according to Lemma \ref{LemInfShiftProps} that $M \subseteq \mathcal X$ is a lattice-infimizer of $f$ if, and only if, $\{0\} \subset C_0^1([a,b]; \R^n)$ is a lattice-infimizer for $\hat f(\cdot; M)$. The following result is a straightforward consequence of Lemma \ref{LemInfimizerScalar} and Proposition \ref{PropTransScalar}.

\begin{corollary}
\label{CorScalarInfimizer}
If $M \subseteq \mathcal{X}$ is a lattice-infimizer of $f$, then $0 \in C_0^1([a,b]; \R^n)$ is a minimizer of $\hat \vp_{f, \zeta}(\cdot; M)$ for every $\zeta \in C^+\bs\{0\}$.

Conversely, if $M \subseteq \mathcal X$ and $0 \in C_0^1([a,b]; \R^n)$ is a minimizer of $\hat\vp_{f, \zeta}(\cdot; M)$ for every $\zeta \in C^+\bs\{0\}$, then $M$ is a lattice-infimizer of $f$.
\end{corollary}

In the following, the derivative (gradient) of the function $\zeta^T L \colon [a,b] \times \R^n \times \R^n \to \R$ with respect to the second and third variable is denoted by 
\[
\frac{\partial}{\partial y} \zeta^T L(t, y, p) \quad \text{and} \quad \frac{\partial}{\partial p} \zeta^T L(t, y, p),
\]
respectively.

\begin{proposition}
\label{PropOC-CVP}
Let $L \colon [a,b] \times \R^n \times \R^n \to \R^d$ be of class $C^1$ and $M \subseteq \mathcal{X}$. Assume that there exists a solution $x_\zeta \in M$ of the problem
\[
\tag{ScCVP} \text{minimize} \quad \vp_{f, \zeta}(x) = \zeta^T F(x) \quad \text{subject to} \quad x \in \mathcal X
\]
for each $\zeta \in C^+\bs\{0\}$. Then $M$ is a sc-solution of (CVP), $\hat \vp_{f, \zeta}(\cdot; M)$ is Fr\'echet-differentiable at $0$ on $C_0^1([a,b];\R^n)$ and
\[
\hat \vp_{f, \zeta}'(0; M)(u) = \int_a^b\left[\frac{\partial}{\partial y} \zeta^T L(t, x_\zeta, \dot x_\zeta)u +
    \frac{\partial}{\partial p} \zeta^TL(t, x_\zeta, \dot x_\zeta) \dot u \right] dt = 0
\]
for every $u \in C_0^1([a,b]; \R^n)$ and $\zeta \in C^+\bs\{0\}$.
\end{proposition}

{\sc Proof.} Corollary \ref{CorVectorScalar} gives that $M$ is a lattice-infimizer for the vectorial calculus of variations problem, hence a sc-solution for (CVP).

Take $\zeta \in C^+\bs\{0\}$. One has $\mathcal X = \{x\} + C_0^1([a,b]; \R^n)$ for each $x\in\mathcal X$, hence $x_\zeta \in M \subseteq \mathcal X$ implies $M\subseteq \cb{x_\zeta}+C_0^1([a,b]; \R^n)$.
  
Since the functional $u\mapsto\psi(u):=\vp_{f, \zeta}(x_\zeta + u)= (\zeta^T F)(x_\zeta + u)$ is Fr\'echet-differentiable on $C^1_0([a,b];\R^n)$ with
\begin{align*}
\psi'(0)(u) &=  \int_a^b\left[\frac{\partial}{\partial y} \zeta^T L(t, x_\zeta, \dot x_\zeta)u +
\frac{\partial}{\partial p} \zeta^TL(t, x_\zeta, \dot x_\zeta) \dot u \right] dt
\end{align*}
(see \cite[\S 0.2, Examples 1, 7 and 8 as well as \S 2.2]{IoffeTichomirov79Book}), one can apply Proposition \ref{PropDiffGen} with $D=C^1_0([a,b];\R^n)$ and the statement follows. 
\pend

Taking Corollary \ref{CorVectorScalar} into account one may observe that the first order optimality condition for $\hat\vp_{f, \zeta}(\cdot; M)$, namely $\hat \vp_{f, \zeta}'(0; M)(u) = 0$, produces the same result as the direct (scalar) optimality condition applied to the solutions $x_\zeta$ of (ScCVP). This is, of course, due to the strong assumption that each scalarized problem has a solution. Again, one can try to solve as many of the scalarized problems as possible and check if the resulting set of solutions already is an infimizer.

\section{A remark on maximization}

The transition from minimization to maximization is slightly more involved for set optimization problems compared to the real-valued case. Instead of the inf-translation, the sup-translation
\[
\check f(x;M)= \sup_{y\in M} f(x+y)
\]
of a function $f \colon X \to W$ with values in a complete lattice $(W, \leq)$ and for a subset $M \subseteq X$ of the linear space $X$ has to be used. Lemma \ref{LemInfShiftProps} remains true with the obvious changes from $\inf$ to $\sup$, and in Lemma \ref{LemInfShiftConvex}, convexity for $f$ has to be replaced by concavity. Semiconvexity and inf-additivity of $(W, +, \cdot, \leq)$ in Lemma \ref{LemConvexImages} has to be replaced by semiconcavity and sup-additivity, respectively, but the result remains the same.

The standard image lattices are $(\F(Z,-C), \subseteq)$ and $(\G(Z,-C), \subseteq)$, and
\[
\psi_{f,z^*}(x) = \sup_{z \in f(x)} z^*(x)
\]
is the appropriate scalarization. The basic fact parallel to Lemma \ref{LemInfimizerScalar} is that if $M$ is an $(\F(Z,-C), \subseteq)$- or a $(\G(Z,-C), \subseteq)$-supremizer of $f \colon X \to \F(Z,-C)$, then $0$ is a maximizer of $\psi_{\check f(\cdot; M), z^*}$ for every $z^*\in C^+\bs\{0\}$.  For the (b) part, $f \colon X \to \G(Z,-C)$ has to be concave.

With these and similar changes, the maximization theory becomes completely symmetric to the minimization theory.

\medskip
{\bf Acknowledgement.} The work of D. Visetti was supported within the project {\em Verification Techniques for Multicriteria Variational Problems} by Free University of Bozen-Bolzano (Grant VerTecMVP).

\end{document}